\documentclass{PiInTheSky}
\usepackage{xcolor}
\usepackage{epigraph}
\usepackage{hyperref}
\usepackage{float}
\begin{document}
\definecolor{myblue}{RGB}{70, 40, 230}
\definecolor{mygreen}{RGB}{100,250,0}
\title{\begin{flushleft}\textcolor{myblue}{Topological Shapes and Their Significance:} \\ \large \textcolor{magenta}{Playing with Loops, Scissors and Glue.}\end{flushleft}}
\author{\href{mailto: kaziaburousan@gmail.com}{\textcolor{red}{K.A.Rousan}}}
\maketitle

\epigraph{A child's first geometrical discoveries are topological.If you ask him to copy a square or a triangle, he draws a closed circle.}{--Jean Piaget}
\thispagestyle{empty}

\section{Introduction}
I am sure while solving any problem or studying something, you have imagined(or may have seen any animation) how shapes transform from one to another. Yes, I know it seems kind of cool. But why I am saying this now?, Well as you will come to know they are understand somehow related to Topology (yeah, if it's not then why the hell I brought it up now).\textbf{We can get those transformation using the concept of Topology (Just with some restrictions)}.
\\
So, What is Topology?, In simple terms:It is the branch of Mathematics that studies \textbf{\textcolor{blue}{pattern and relative position, without regrad to size}} (Yes, Topologists believe in equality , so for them size doesn't matter). Topology is sometimes referred to as "\textbf{Rubber-Sheet Geometry}" because a figure can be changed into equivalent figure by \textit{\textbf{bending}, \textbf{stretching}, \textbf{twisting}, \textbf{squishing}, and the like, } but \textit{not by \textbf{tearing} or \textbf{cutting}}(as we don't like violence).

\section{Different Topological Shapes}
Topological shapes are very fun to visualize and when we give them colour, it becomes more beautiful!!. But why should we read this,just for fun?(Come on,most of people don't study mathematics just for fun!!). Then why study it?, Well, it helps us solve different problems (What? shapes help us solving problems?). Yes sort of , I mean Analytical Topology help us solve problems but some times shapes also do the same. But how ?, Well,in this article we will see that for 3 shapes(not how they help, but their correspondence with plane) \textbf{\textcolor{red}{A new shape}, \textcolor{blue}{Torus}} and \textbf{\textcolor{mygreen}{Mobius Strip}.} The cause of choosing this 3 shapes are that they are well known and are simpler to visualize than other shapes(like \textbf{Klein Bottle} or \textbf{Projective Plane}, you know I am also lazy).
\begin{figure}[H]
\includegraphics[width=1.8in]{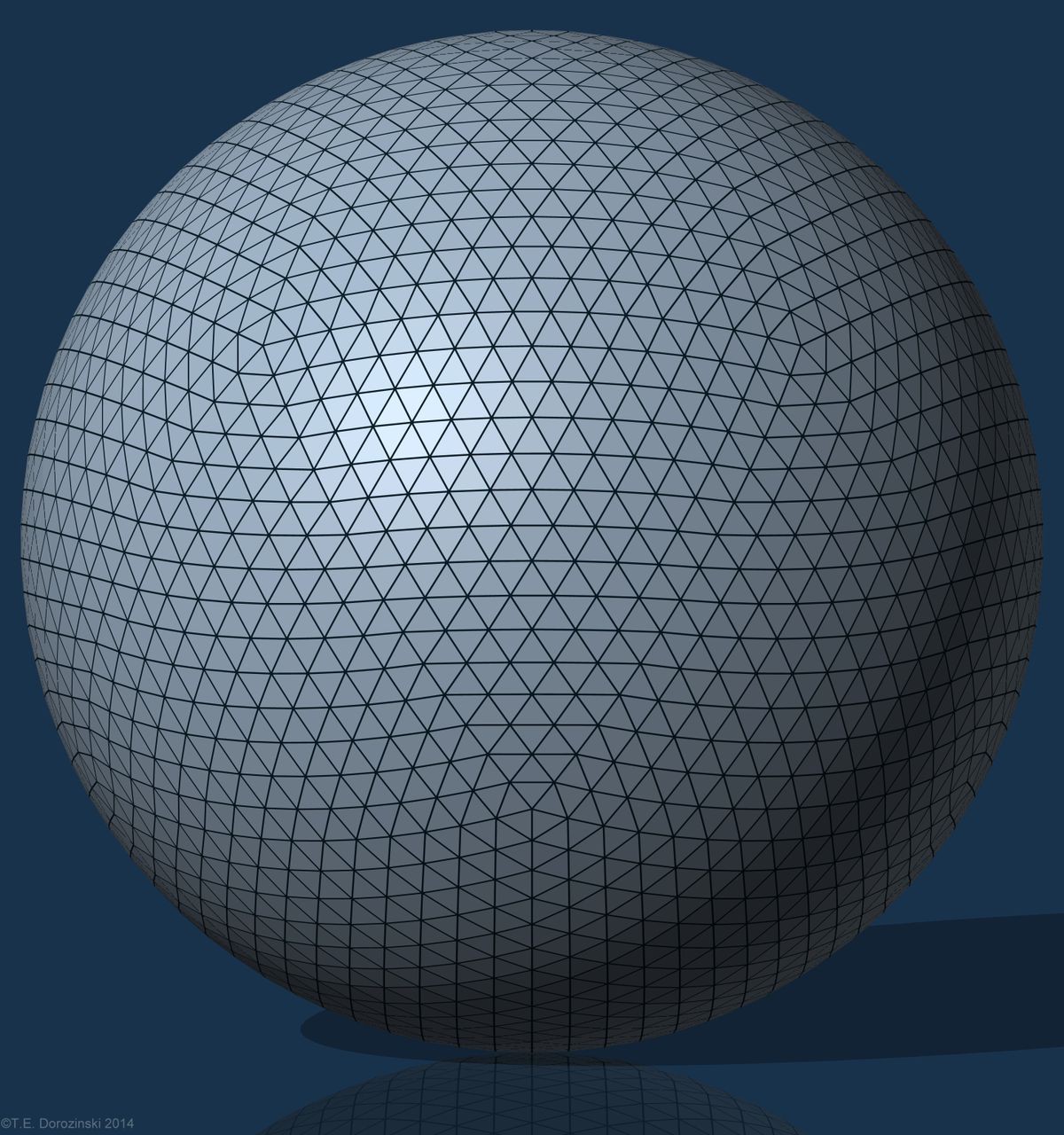}
\caption{Our loving Spherious(we will need this )}
\label{1}
\end{figure}

\begin{figure}[H]
\includegraphics[width=2in]{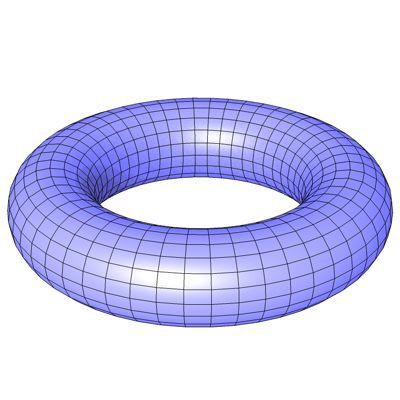}
\caption{Tasty Donut for Topologists}
\label{2}
\end{figure}

\begin{figure}[H]
\includegraphics[width=2in]{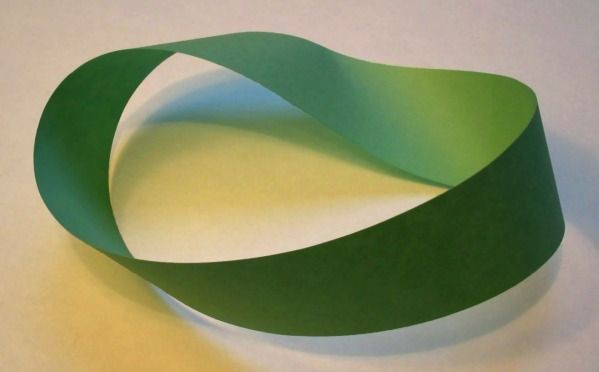}
\caption{Our Mobius(not Ultraman Mibus)}
\label{3}
\end{figure}
Now get ready with a loop, a scissor, some tape and ink. Let's see the beauty of Mathematics\textbf{(}before reading the next section I suggest to read \textbf{The $\pi$ In The Sky Magazine Issu: 21 \textcolor{red}{\href{http://www.pims.math.ca/resources/publications/pi-sky}{Mathematical Cut and Paste, An Introduction to the Topology of Surface}} by Maia Averett}). There you will learn about Topological shapes and how to create them from a 2D shape , which you will also see here,but if you have read that then it will help you a lot to understand to understand whatever we will discuss here.
\section{Shapes and Hidden-maths}
\begin{figure}[H]
\includegraphics[width=2in]{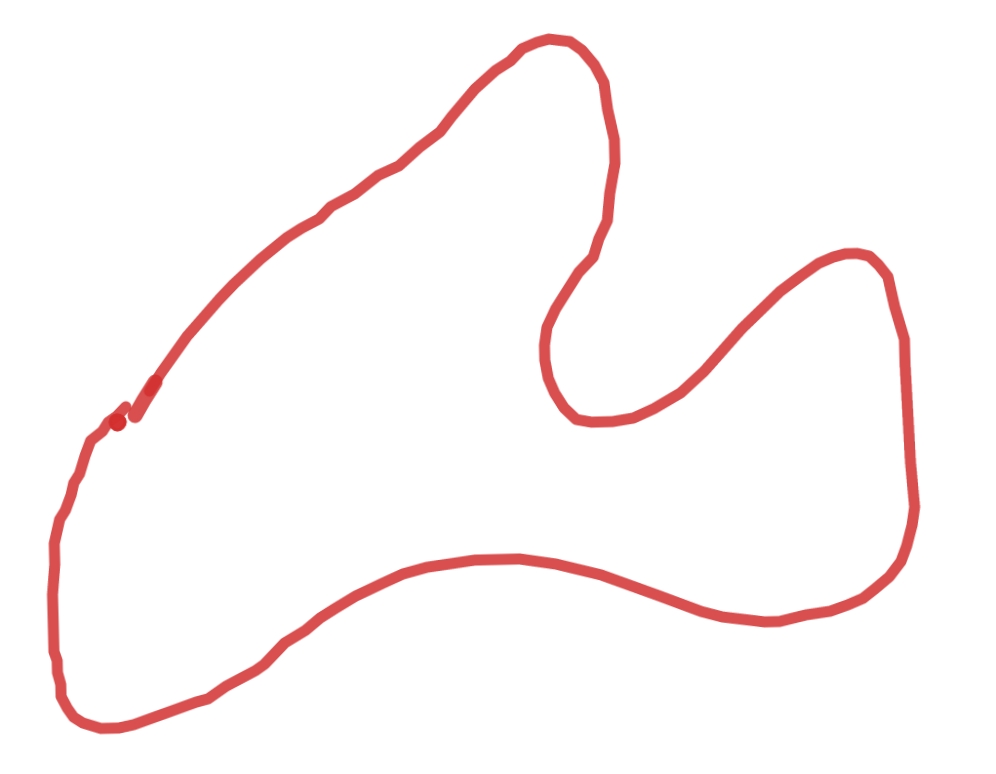}
\caption{Any random loop}
\label{loop}
\end{figure}

So what do you think,how can we get this shapes make us help?(We shall see here only one form of there hospitality on how they save us from different problems).  \\
Suppose you have a loop; any random 2D loop is fine as shown in fig-\ref{loop}(don't laugh, On my drawing skills).
\paragraph{Hunting a funny shape:} Now suppose we have a \textbf{point A} outside the loop and a \textbf{point B} inside the loop and we seek the information regarding those two points. Now notice something important; we can make a ordered pair like this \textbf{\textcolor{blue}{an pair like this (Position of A, Position of B)}(You may ask an ordered or unordered?, Here our main concern is about ordered pair i.e., for any pair $(a,b)\neq (b,a)$)}. So how nice it would be if we can somehow get this information of $(A_p,B_p)$ (where $A_p$ and $B_p$ represent the position of A and B) with just a single surface, with keeping all the continuous things continuous (we math lovers are obsessed with continuous things). \textbf{Like, suppose we make a surface on which each point correspond to any possible ordered pair of numbers on the loop.}
\begin{figure}[H]
\includegraphics[width=2in]{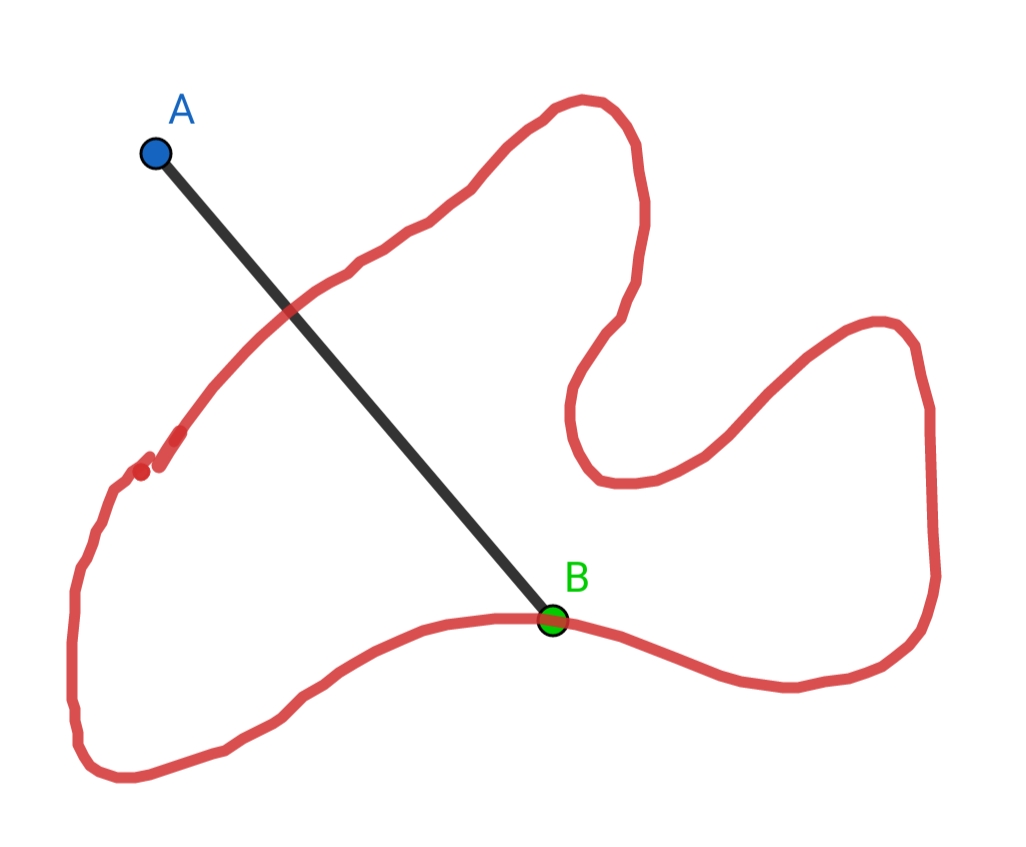}
\caption{The loop with B on It.}
\label{loop1}
\end{figure}

\begin{figure}[H]
\includegraphics[width=3in]{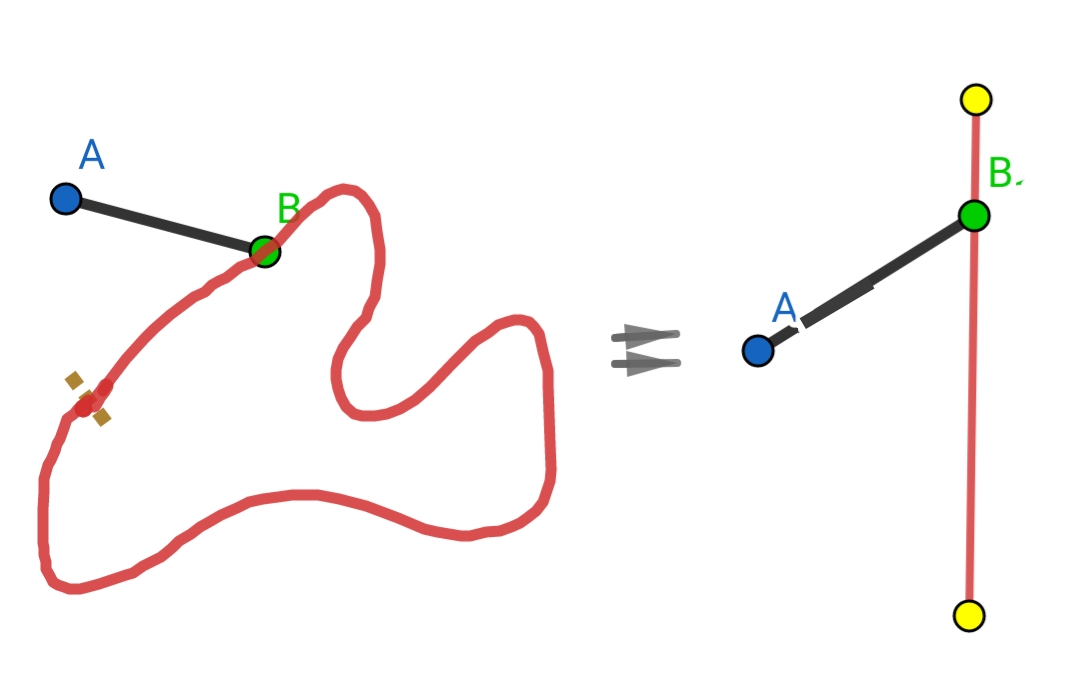}
\caption{The loop is cut along the brown dotted line and The yellow points on the straight line are one single point.}
\label{loopcut}
\end{figure}

\begin{figure}[H]
\includegraphics[width=2.3in]{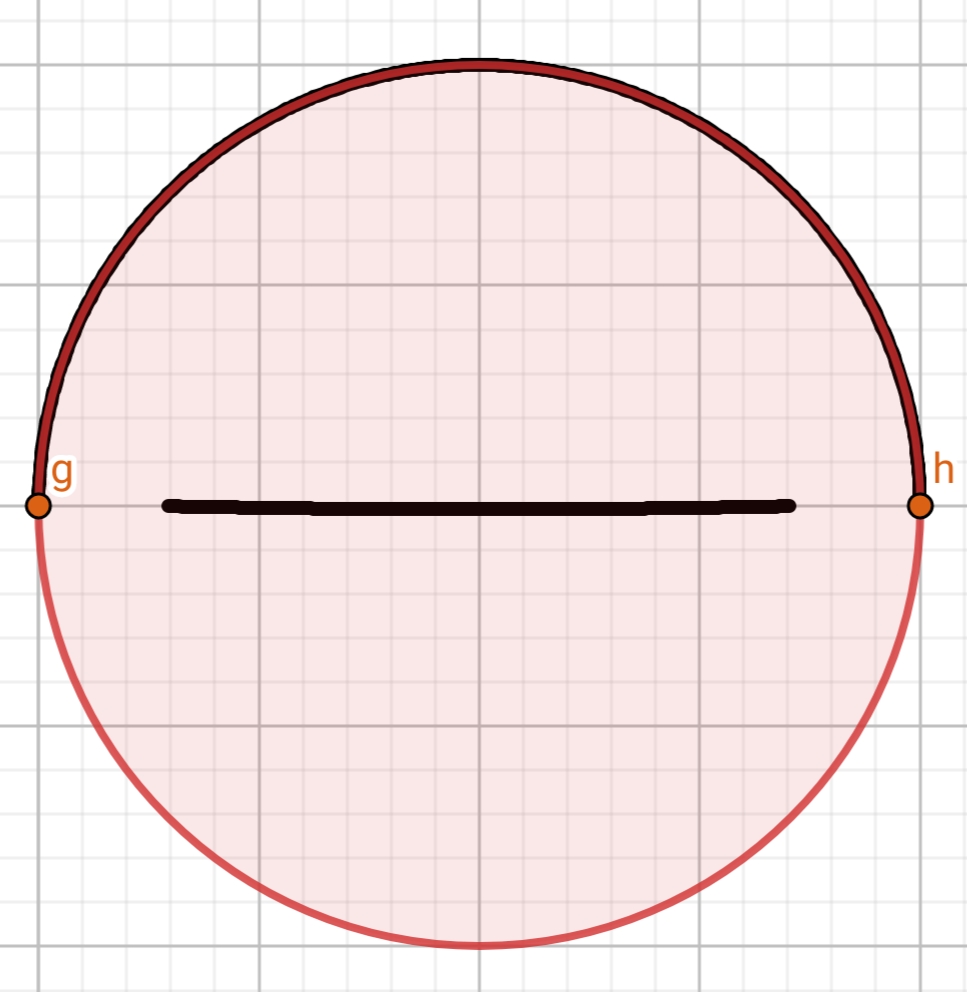}
\caption{It is what we get after we condense those lines.}
\label{circlesp}
\end{figure}

As shown in Fig-\ref{loop1} we have the points where \textit{A is a fixed one but B can move all over the loop}. Now we can do something cool , we can cut any point of the loop and make the loop straight. Doing such we can easily manipulate the loop but as we don't like violence,so we have to \textbf{glue it back such that the point where we cut the loop correspond to the same place}.\\
Here you have to remember one thing: \textbf{The terminal points of the line, which seems to be different, are actually the same point(as we had cut a single point and make 2 terminals)}(see Fig-\ref{loopcut})
Now we can stretch out the point A as a line of equal length of that of the straightened loop. Using this as \textbf{x} and \textbf{y} axis we make a plane, where each point corresponds to each \textbf{possible pair of point of A and B}. 
\begin{figure}[H]
\includegraphics[width=2.5in]{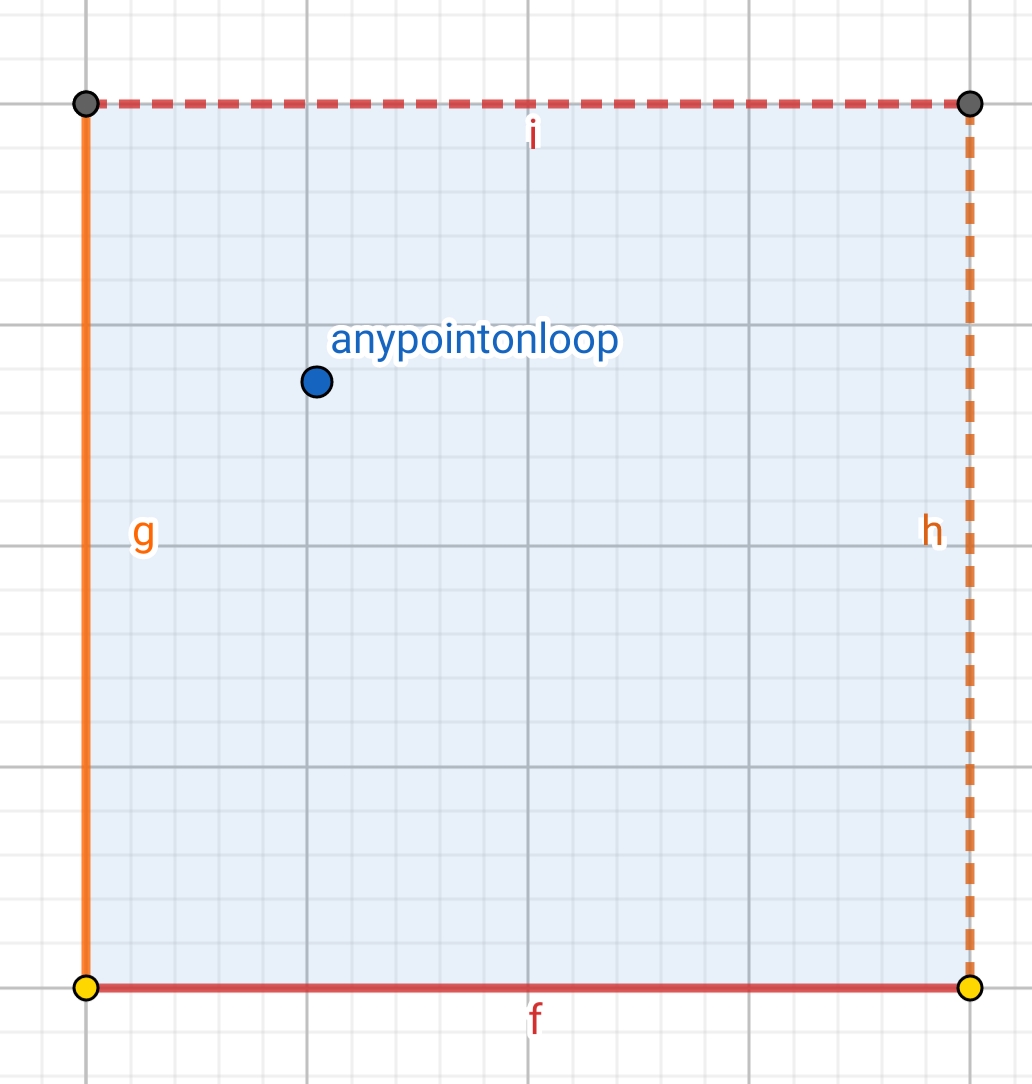}
\caption{We get this plane area from our argument.}
\label{square1}
\end{figure}
In Fig-\ref{square1} the g and h are equivalent line and so are f and i, i.e., they $g\equiv h$ and $i\equiv f$. So \textbf{at the end we have to glue g and h and also i and f in the right alignment.} Also remember g and h are actually a point so the whole line will be condense to a single point and we shall be adding those two points together(this also save us one more work as we also have to add the 2 yellow point, remember?). \textcolor{blue}{So condensing g and h give us a circle.} Now the upper half(upper part of the black line) should be added to it's below half(see Fig-\ref{circlesp}).

If we add those two parts then we will get a beautiful Sphere. Our loving Sphere.
\begin{figure}[H]
\includegraphics[width=2.4in]{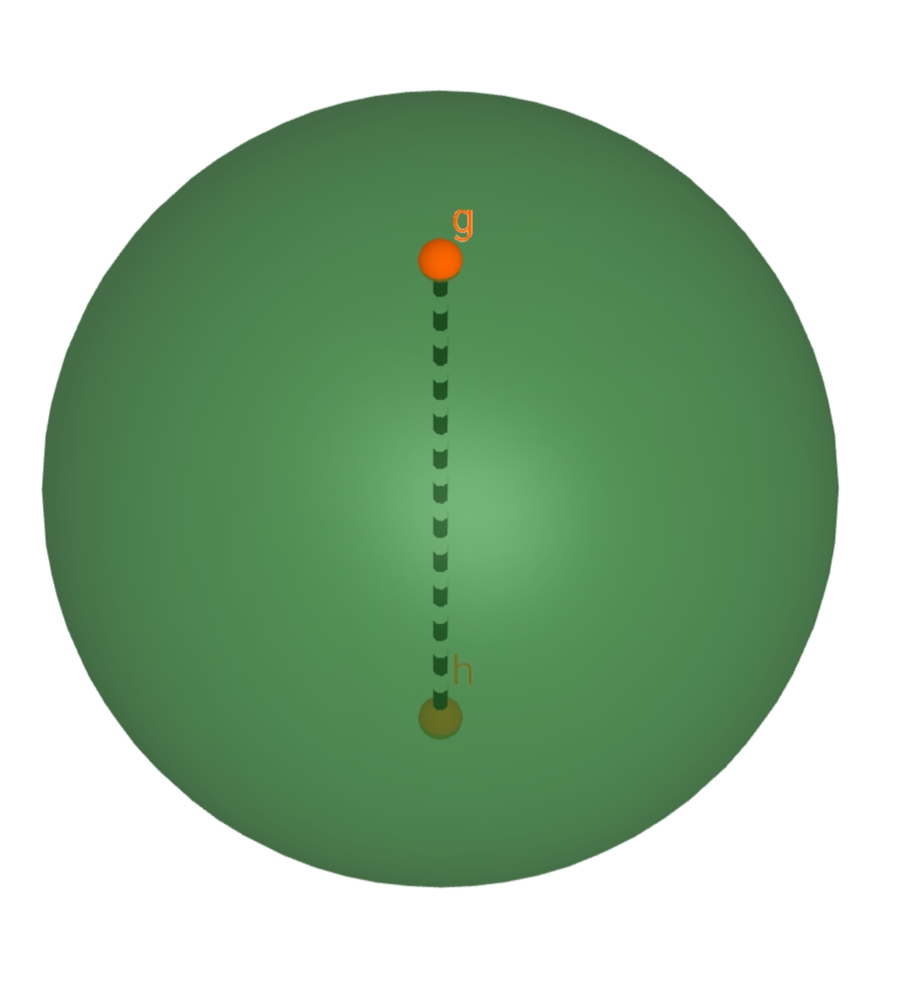}
\caption{This is what we get and h and g are diametrically opposite.}
\label{sphere1}
\end{figure}
Now remember we have to join \textbf{g} and \textbf{h} as they are the same point. So imagine, we have to put our fingers on those points and then press them so that they touch eachother and then freeze it there. The shape we get showhow look like this.

\begin{figure}[H]
\includegraphics[width=2.5in]{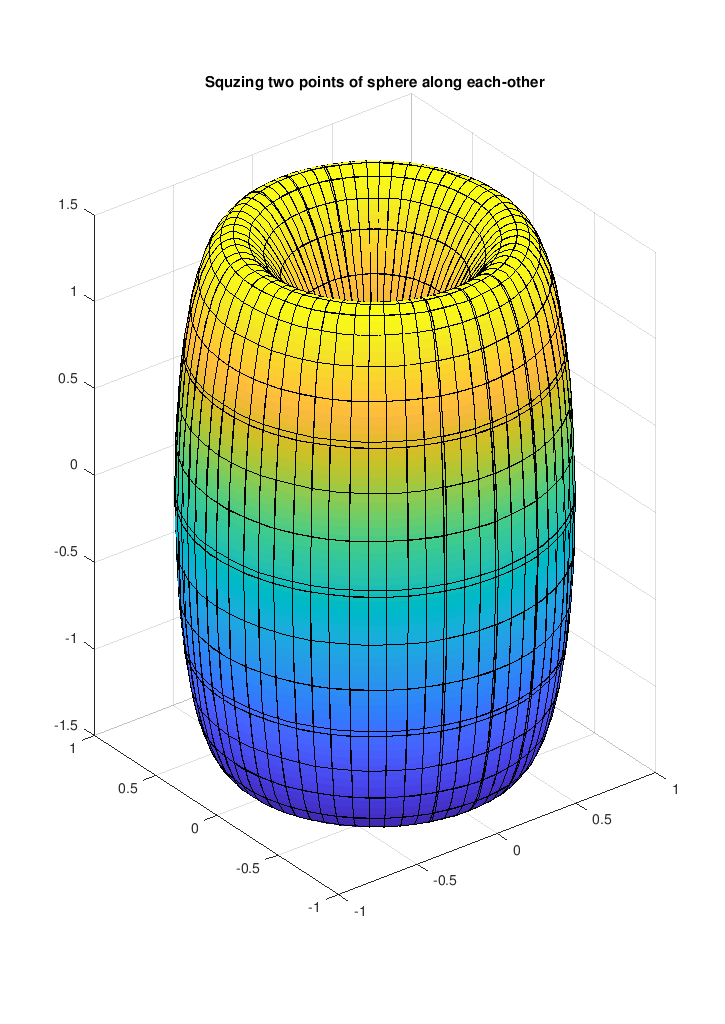}
\caption{This is what we get when we press the points g and h. We then glue those 2 points together }
\label{ST}
\end{figure}
\begin{figure}[H]
\includegraphics[width=2.7in]{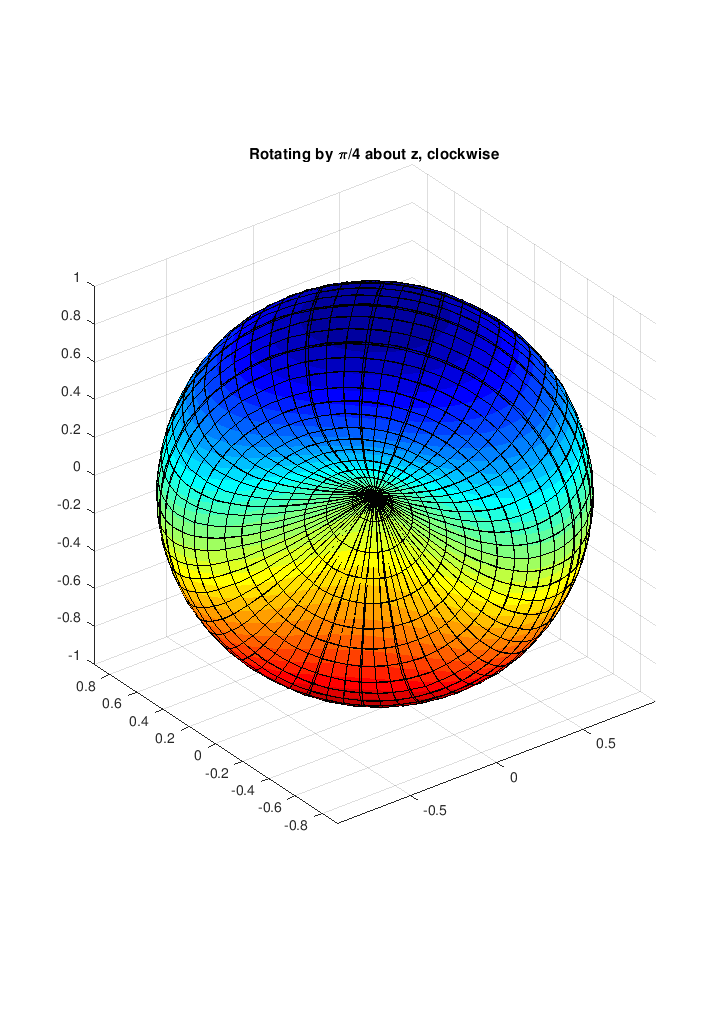}
\caption{Adjusting Fig-\ref{ST} and then rotating it , see it doesn't have any hole.}
\label{ST2}
\end{figure}
\begin{figure}[H]
\includegraphics[width=2.6in]{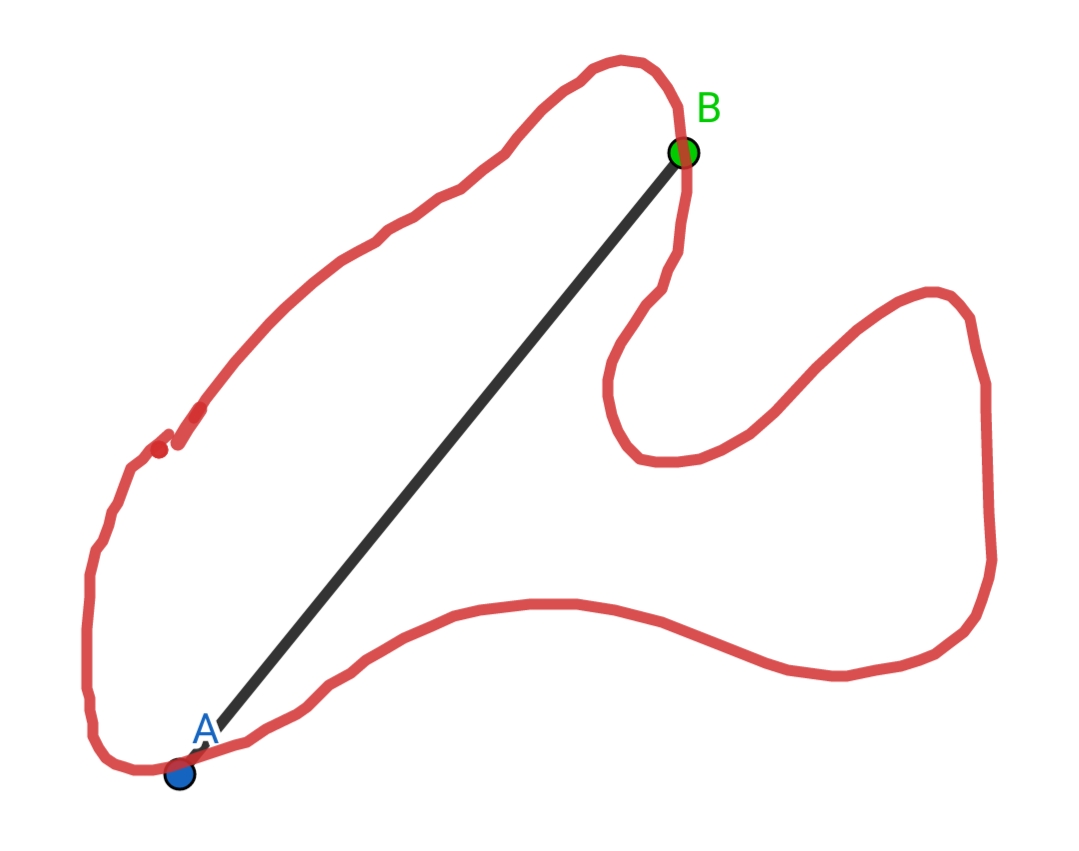}
\caption{Here both points are on the loop}
\label{Tor}
\end{figure}
\begin{figure}[H]
\includegraphics[width=2.5in]{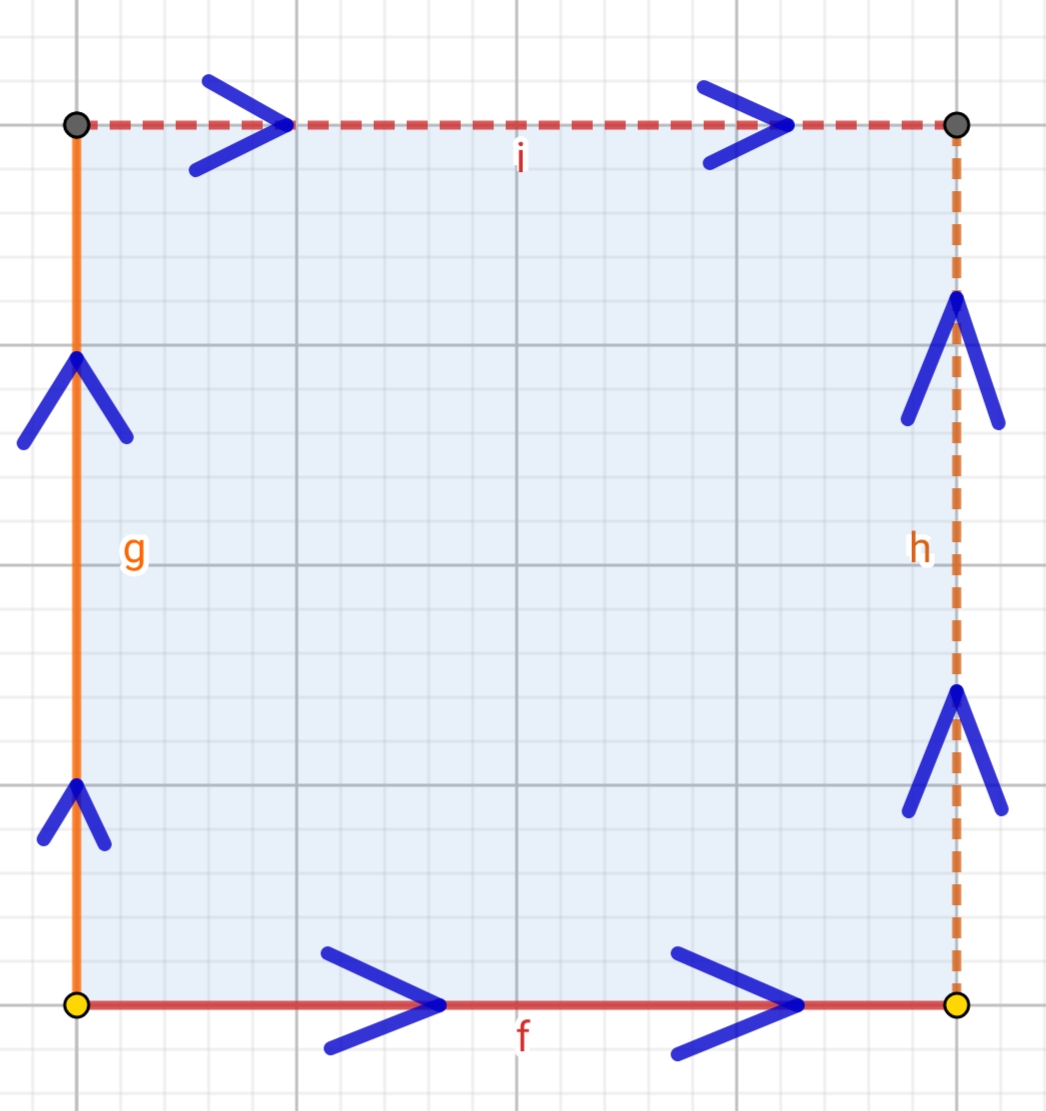}
\caption{Here we have to add the lower f to upper i and right g to left h.}
\label{Torpla}
\end{figure}
\begin{figure}[H]
\includegraphics[width=3in]{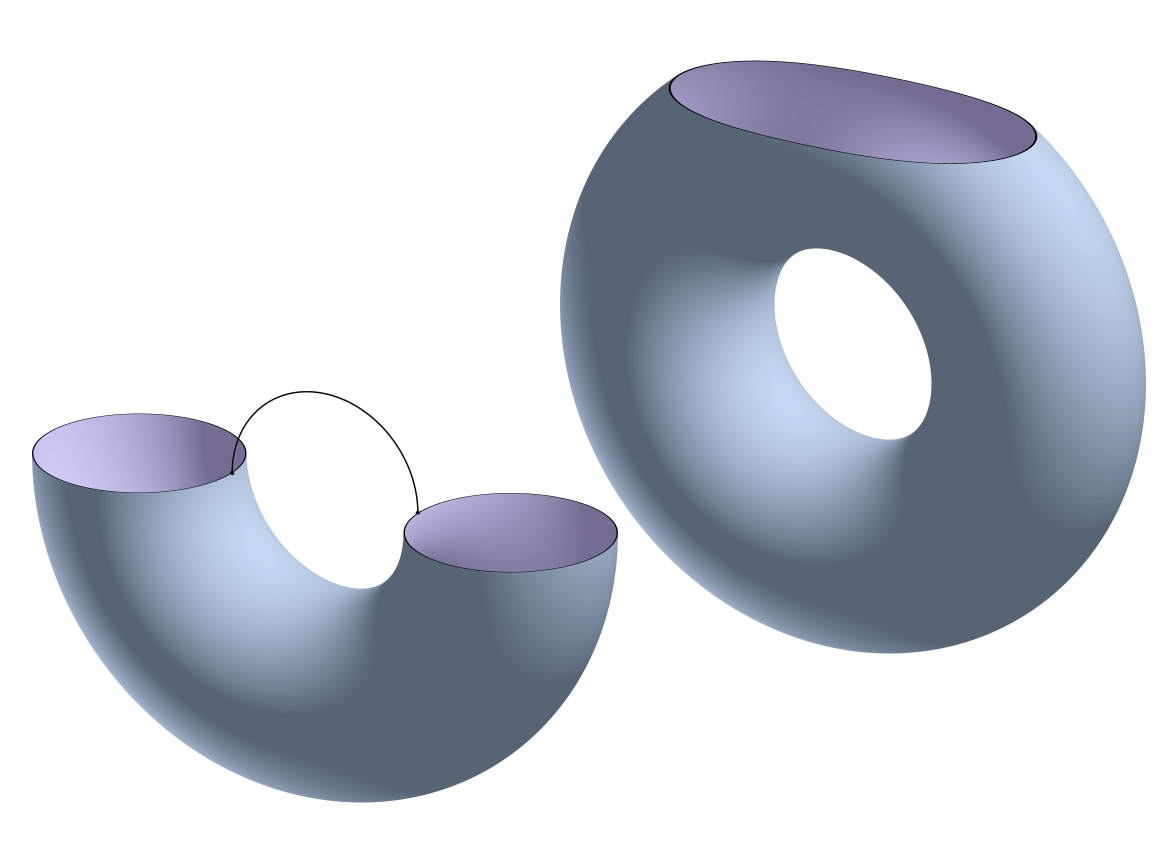}
\caption{Here you see how cylinder curves and slowly become Torus.}
\label{cytor}
\end{figure}

Now recall we just care about shape and we can deform them freely!!, That means we can now deform this shape in figure-\ref{ST} and we will get something like Torus(Fig-\ref{2}). \textcolor{blue}{But it is not the same as \textbf{Torus (which has a hole in it).As This shape has not any hole, it is just squized as shown in }fig-\ref{ST2}}.
\textbf{What do we get from all this?}\\
If we take any point on the this surface (Fig-\ref{ST2}) that point corresponds to the ordered pair of The point on the curve(Fig-\ref{loop}) and outside the curve. So we don't have to care about the points(2point) anymore.All we have to do is to see if our point(one single point!!) is on the surface that we get. But why should we care about that? Well there are 2 reasons.\\
1.It makes for work easier.\\
2.It helps us in certain problems.
\paragraph{Torus and it's Surface:}
We again take the loop in Fig-\ref{loop}. \textbf{This time we have A and B both on it's surface.}Now as before we want a surface which will give us a direct correspondence with the two points(as in Fig-\ref{Tor}).
Now like before, We shall cut this loop at some point and straight it. As both points are on the same line,so to make a 2D plane we have to making it's copy(The arrows show the direction along which we have to glue later).
\begin{figure}[H]
\includegraphics[width=2.8in]{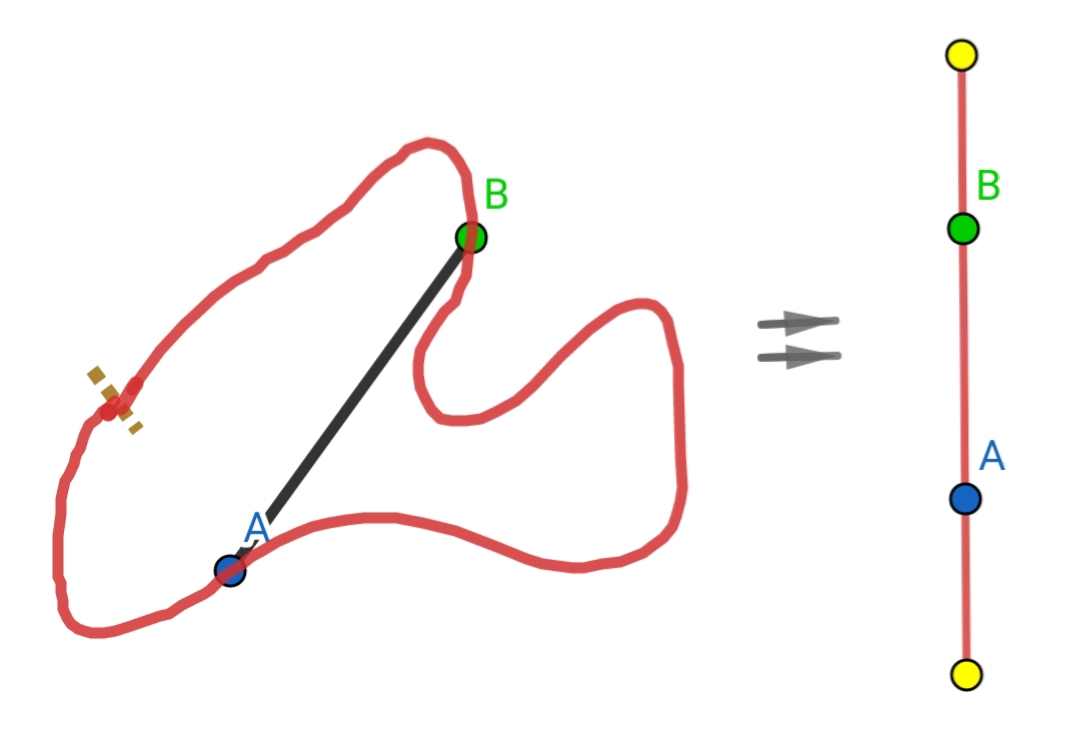}
\caption{See here we have a single line and both points are on it.}
\label{Torloop}
\end{figure}
Here we will use ordered pair(Yes again!). Now making a copy of the line in Fig-\ref{Torloop} and \textbf{taking A on one of the line and B on another} we get a whole new version of our figure.
\begin{figure}[H]
\includegraphics[width=2.5in]{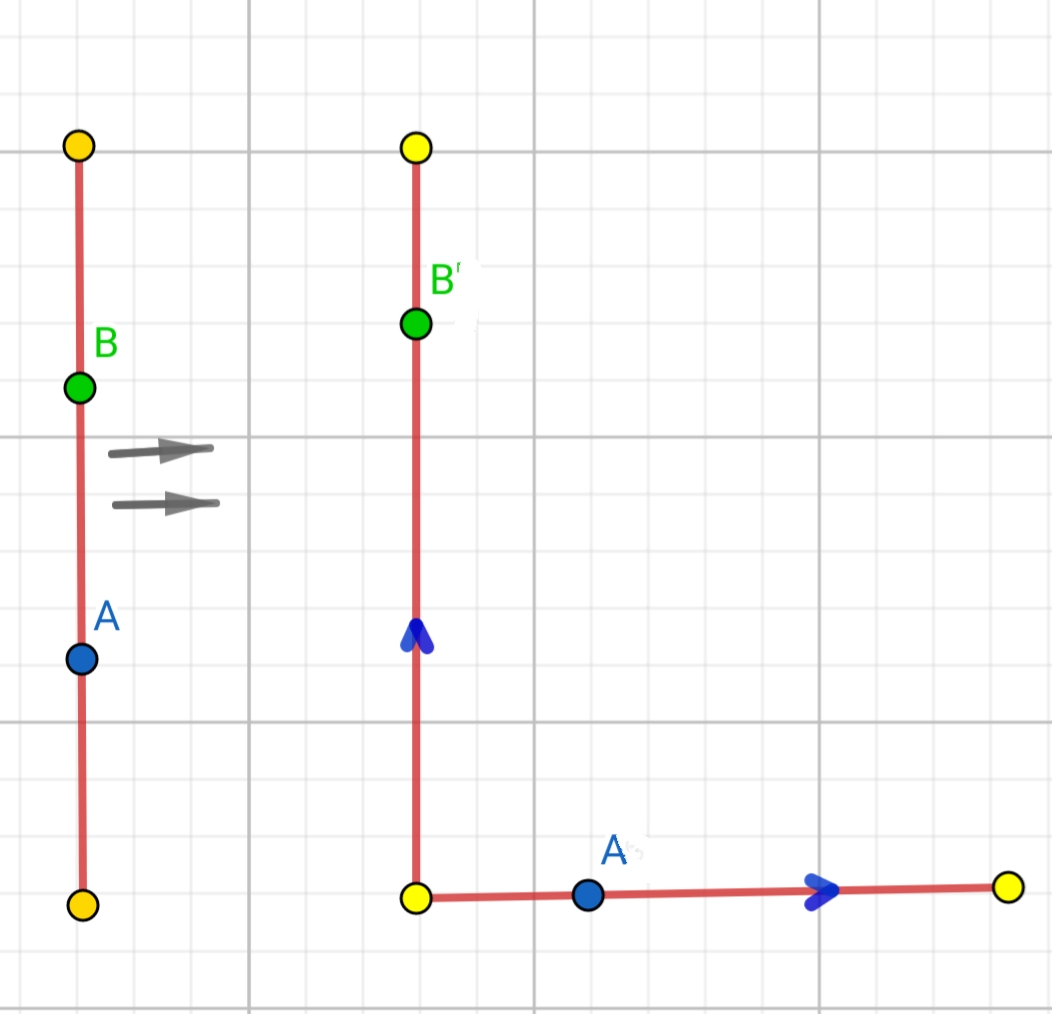}
\caption{Here we have the line looking like this and We can create a plane now using this 2 line as axis.}
\label{Torgrid}
\end{figure}
Now we copy the Horizontal and Vertical lines. Then place those copied ones parallel to their original ones, getting $f\equiv i$ and $g\equiv h$. But we have to paste them in right order as guided by the arrows in Fig-\ref{Torpla}.
Adding the f and i(red lines) we get a cylinder and joining the both ends of the cylinder we get Torus(as in Fig-\ref{cytor}...wow! what a beautiful sight it is). Now if \textbf{We choose a point on the Torus and get the two points on our loop.} Mathematics really is inconceivable.
\paragraph{Mobius Strip's Surface:} 
So for pervious cases all we did is just examine \textbf{Ordered pair}. Now we shall examine unordered pairs(Yes, as if we give too much care on just one other one will become angry!!). Like before we take the loop in Fig-\ref{loop} and we take two points on the loop as in Fig-\ref{Tor}. Now like we did before, Cut the loop at a point and make it straight as in Fig-\ref{Torloop}. Now make two mutually perpendicular axis with it and then make a plane with appropriate arrows(Fig-\ref{Torgrid} and Fig-\ref{Torpla}).
\\
I am sure you are thinking what a lazy Author, He is just giving us directions. yes I am as \textbf{the first few steps are similar.}\\
Now as I said before We shall be using Unordered pair. So for this we have to \textcolor{blue}{keep this in mind that from now on $(2,3)=(3,2)$ i.e., $(x,y)=(y,x)$. But if that's the case then, We have the magical line $x=y$(It's magical as for ordered pair this is the line of symmetry)}. Let's Draw that the line.
\begin{figure}[H]
\includegraphics[width=2.6in]{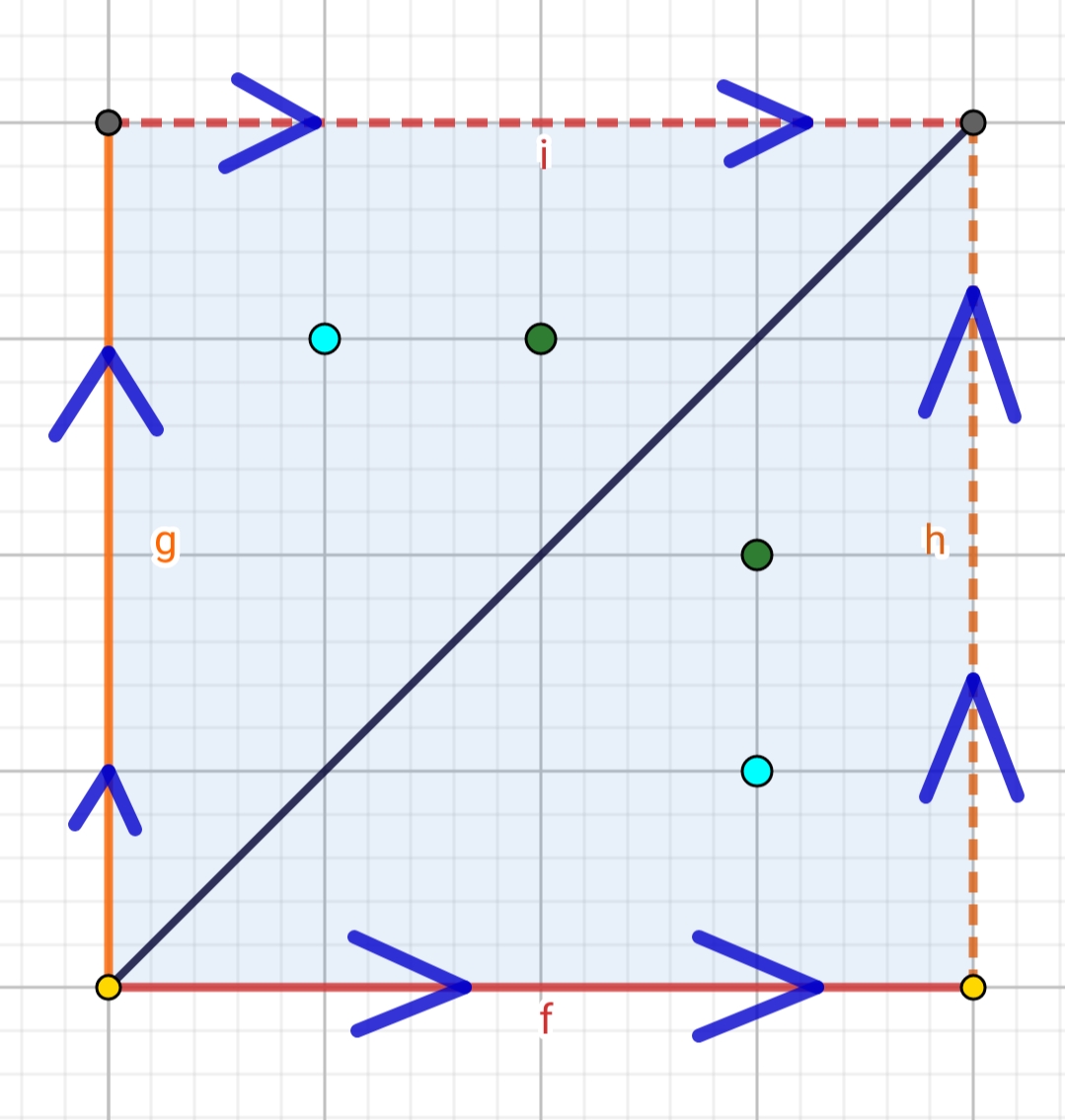}
\caption{Notice that the Unordered pairs are symmetric with respect to $x=y$ line.}
\label{unorder}
\end{figure}
So that means \textcolor{blue}{\textbf{If we choose $(3,4)$ and then choose $(4,3)$,we are choosing same point twice. So it's not unique anymore}.} To keep things unique and not to take same point twice, we have to remove one of the points or fall $(a,b)$ to $(b,a)$ for all \textbf{a} and \textbf{b} then the problem can be solved!!\\
But how??, Well remember I have told you that those points are symmetric about $x=y$. So we can fold the plane along the black line (Fig-\ref{unorder}) and our problem is solved. But now how can we align the arrows while keeping their direction right?,Our only option is to tear things. But we can't just tear things(as it will make things discontinuous and in somecases divergent),so \textbf{\textcolor{blue}{If we tear things we have to paste them in right places in such a way that all things stay continuous and nice.}}

\begin{figure}[H]
\includegraphics[width=2.3in]{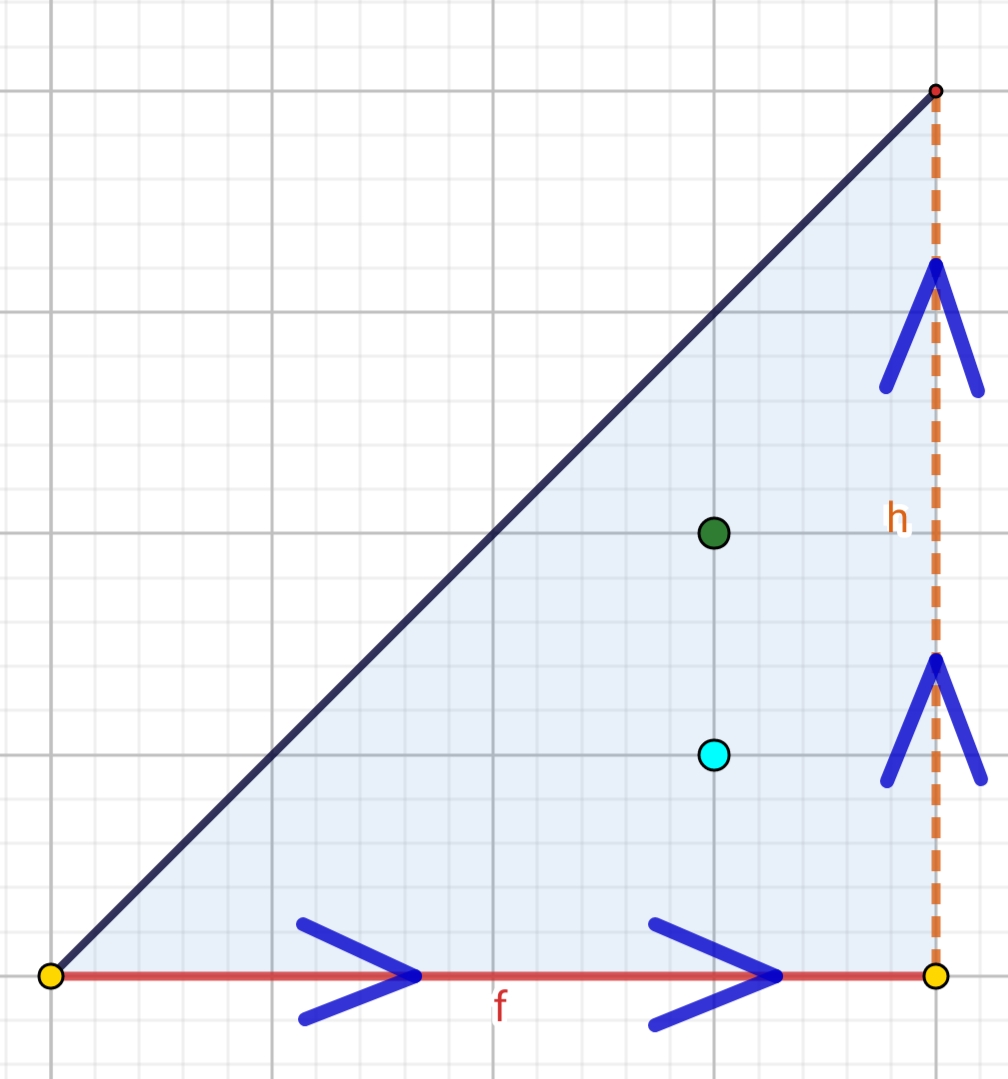}
\caption{Folding two parts along $x=y$ line.}
\label{fold}
\end{figure}

\begin{figure}[H]
\includegraphics[width=2.5in]{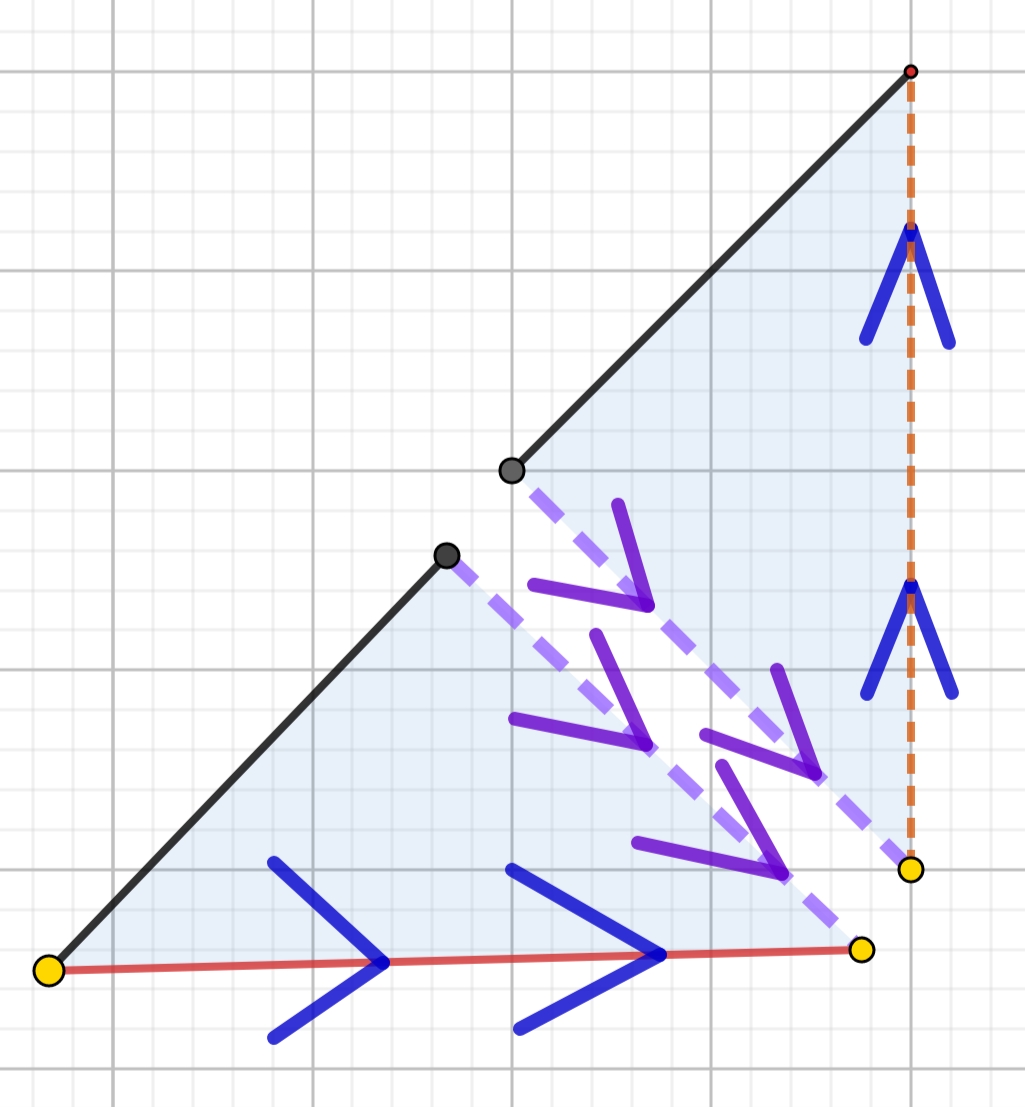}
\caption{The guide line for out cut.}
\label{Cutted}
\end{figure}

\begin{figure}[H]
\includegraphics[width=2.6in]{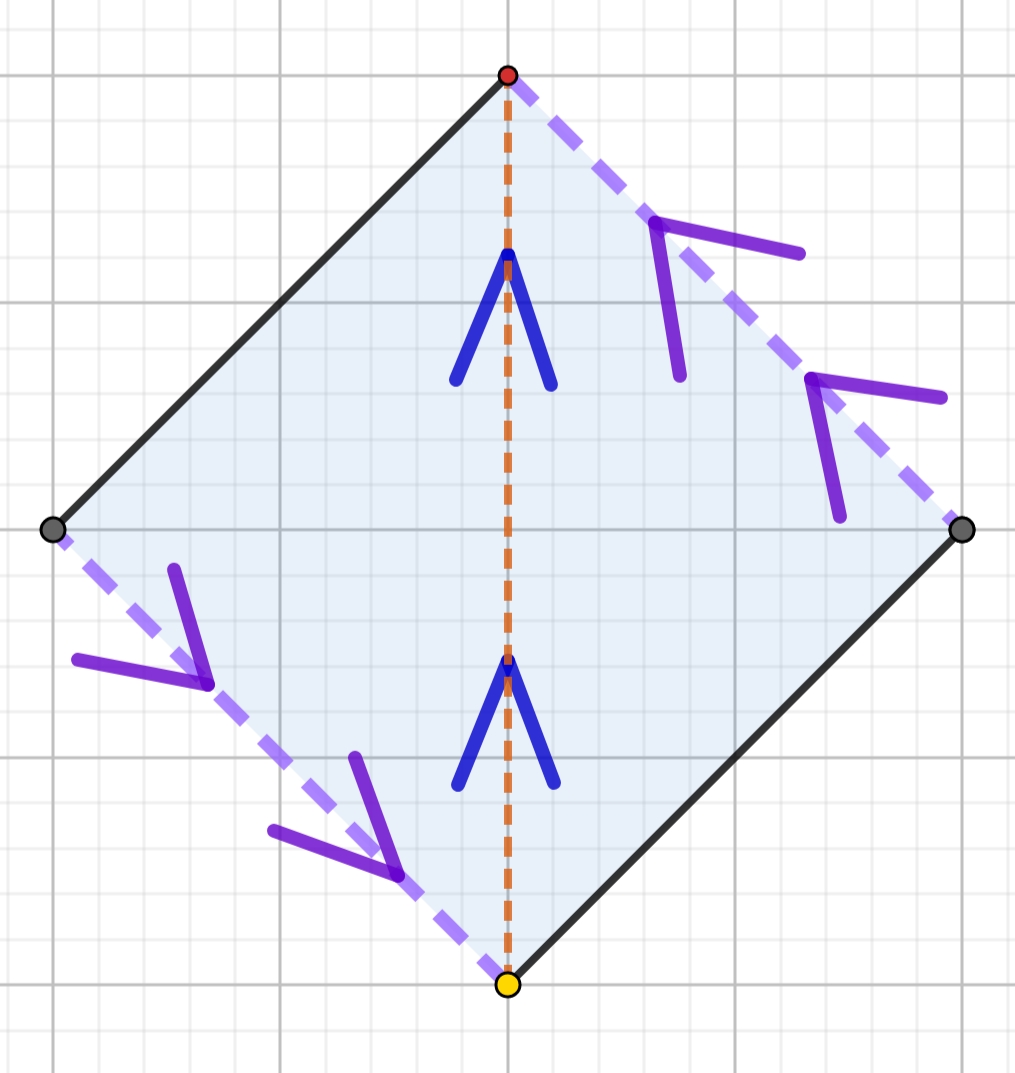}
\caption{The red and organe line are joined along blue arrows.}
\label{Cutjoined}
\end{figure}

Now let's join them according to our guiding arrows. But how can we cut those?, First we have to make a cut as without further dividing, our ultimate goal can't be reached. So let's go it.
\begin{figure}[H]
\includegraphics[width=2.5in]{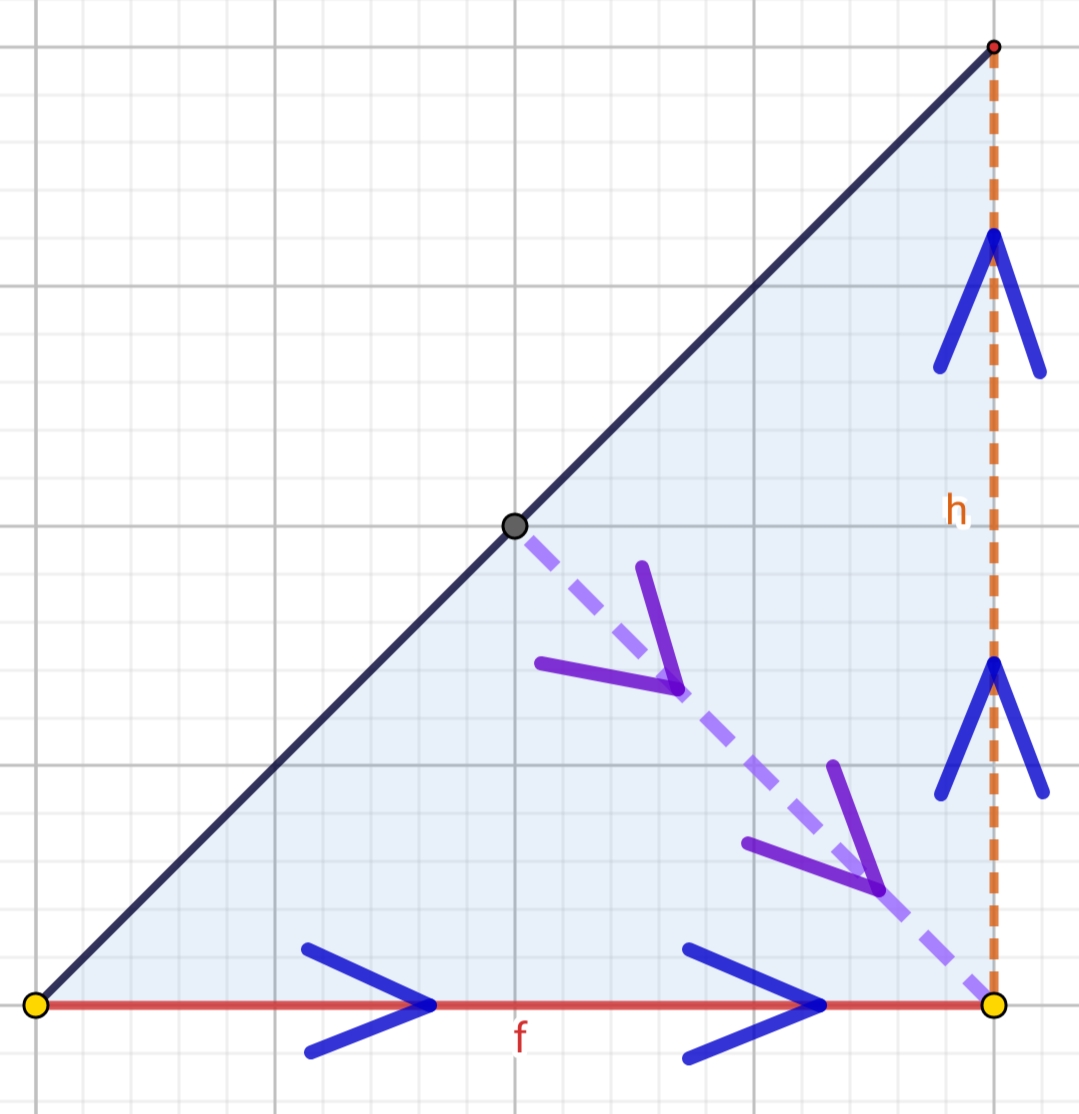}
\caption{The guide line for out cut.}
\label{cutline}
\end{figure}
Now as shown in the Fig-\ref{cutline} the \textbf{violet line} is our line along which we will cut the \textcolor{blue}{\textbf{folded} plane}.The arrows show us in which direction we have to glue things together. After we have completed our cutting, we will get something like Fig-\ref{Cutted}.
Now after judging by the Arrows we see that \textbf{\textcolor{blue}{we have to add the red line with the dotted orange line(The lines with blue arrows)}.}Here one things should be noted that \textbf{the black line represent the $x=y$ line}, i.e., The line where elements like $(a,a)$ exists.\\
Now let's join along the line. While after joining along the line we get a plane as shown in Fig-\ref{Cutjoined}.

\begin{figure}[H]
\includegraphics[width=2in]{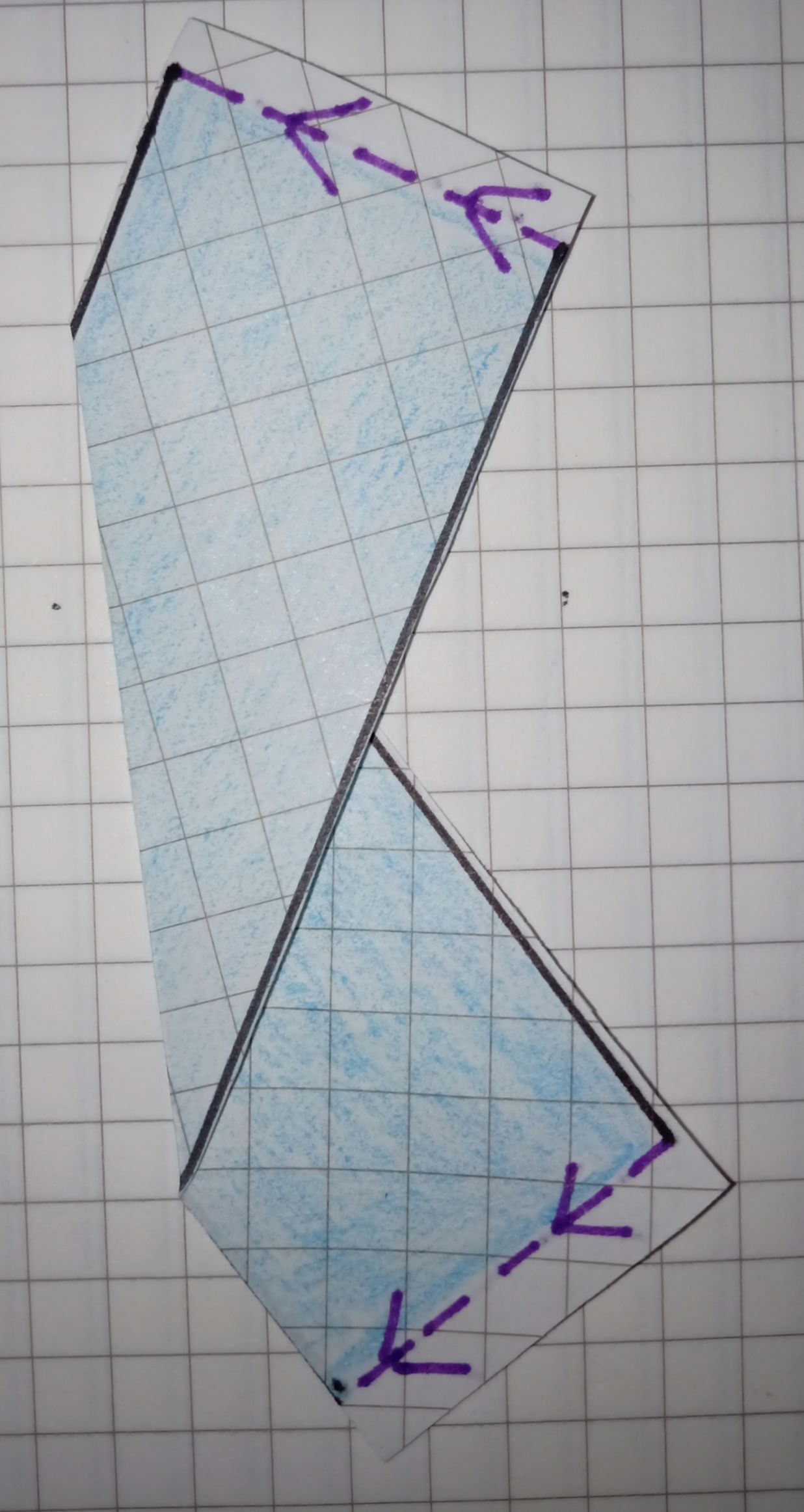}
\caption{One side is rotated $90^{\circ}$ and due to small size of the paper the uneasy fold}
\label{Rotated}
\end{figure}

\begin{figure}[H]
\includegraphics[width=2.5in]{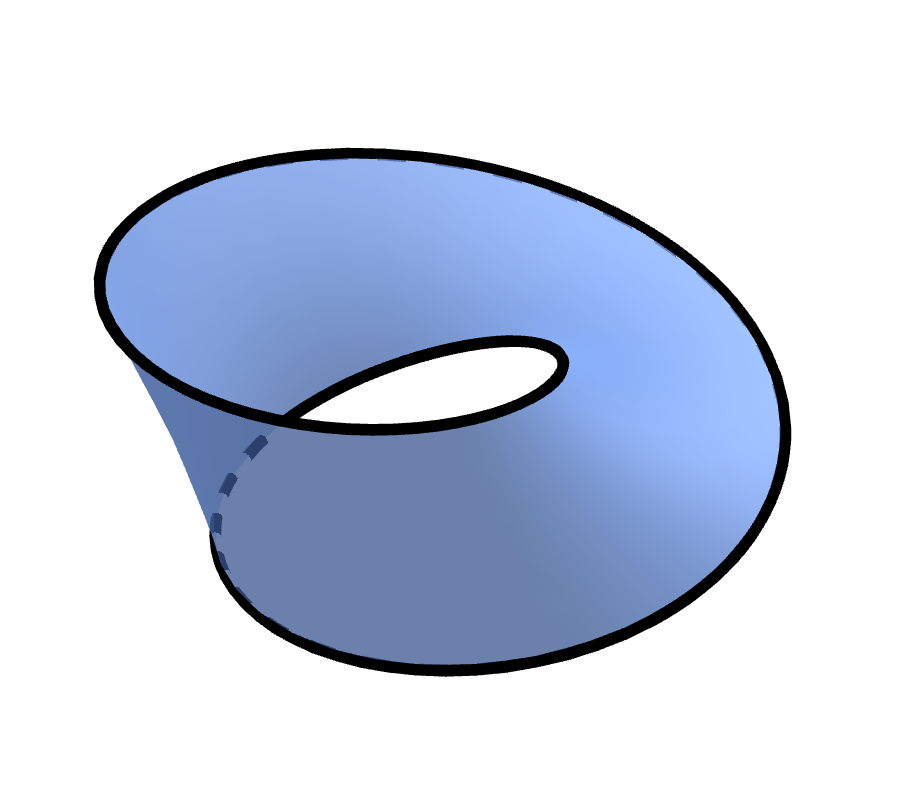}
\caption{Joining the violet line we get this , Mobius Strip}
\label{Final}
\end{figure}

Now we have another pair of lines to join(The violet ones). But we can't join them directly because of the direction of arrows. Then what!!!\\
Well it's simple. First stretch it and then \textbf{Just turn one side with the violet colour $90^{\circ}$}(as shown in Fig-\ref{Rotated}, It is hand drawn as for someone's request and as she used small paper so the rotation is not done clearly), Wow.. The direction of the arrows are similar to the arrows of the other end.
Now after this just join the violet sides(Fig-\ref{Final}). Here what we get is the \textbf{Mobius Strip (whose edge is the line $x=y$, that black line!!)}.
So this is the surface for unordered pair of numbers on the loop. On this surface take any point and that represent two points on our loop in Fig-\ref{loop}.
All the work that we did, are for what?, As I have told before, All those things which we did for all this pages are not just for minimizing our notation but to solve different problems. Here are some problems which can be solved in many ways, but if we use the relationship of Topology to the unordered pair and ordered pair the solution of these problems become super awesome. \\
1. The Inscribed Square Problem.\\
2. The Necklace Problem.\\
3. The Inscribed Triangle problem.\\
4. A problem about parallel Lines...etc,\\
\textbf{Those solutions of this problems can truly be called one of the most elegant solution if we use the significance of 3D Surfaces.} In this article we will not discuss about those problems but I am sure those solutions will surely surprise you.\\
I hope you all liked this article and from this you may have learn something new and wonderful. If you want to learn more go through the references.
\section{Summary}
So From what we learn we can say If we have any Paper then by gluing them in right way, we can create many Topological Surfaces. These surfaces have special \textbf{\textcolor{blue}{Topological Properties which are invariant understand stretching, twisting,.. etc(without tearing).}} When gluing we have to consider the sense of direction. Using this shapes we can define ordered and unordered pairs and using this concept we can \textit{\textbf{assign single point on surfaces to different points}(two or more)}. This helps us in solving many problems which are almost impossible to solve. Here we have seen pressed Sphere correspond to ordered pair for a curve.(a fixed and a variable point on the curve), Torus also corresponds to ordered pair(here both points can move along a single curve) and finally Mobius Strip is related to unordered pair(Same as Torus).

\section{Acknowledgement}
I would like to thank friends Syeda Spandita Zaman and Bipul Karmaka and also to Pinaki Ranjan Ghosh to help me writting this article.

\picturebiography{Kazi Abu Rousan}{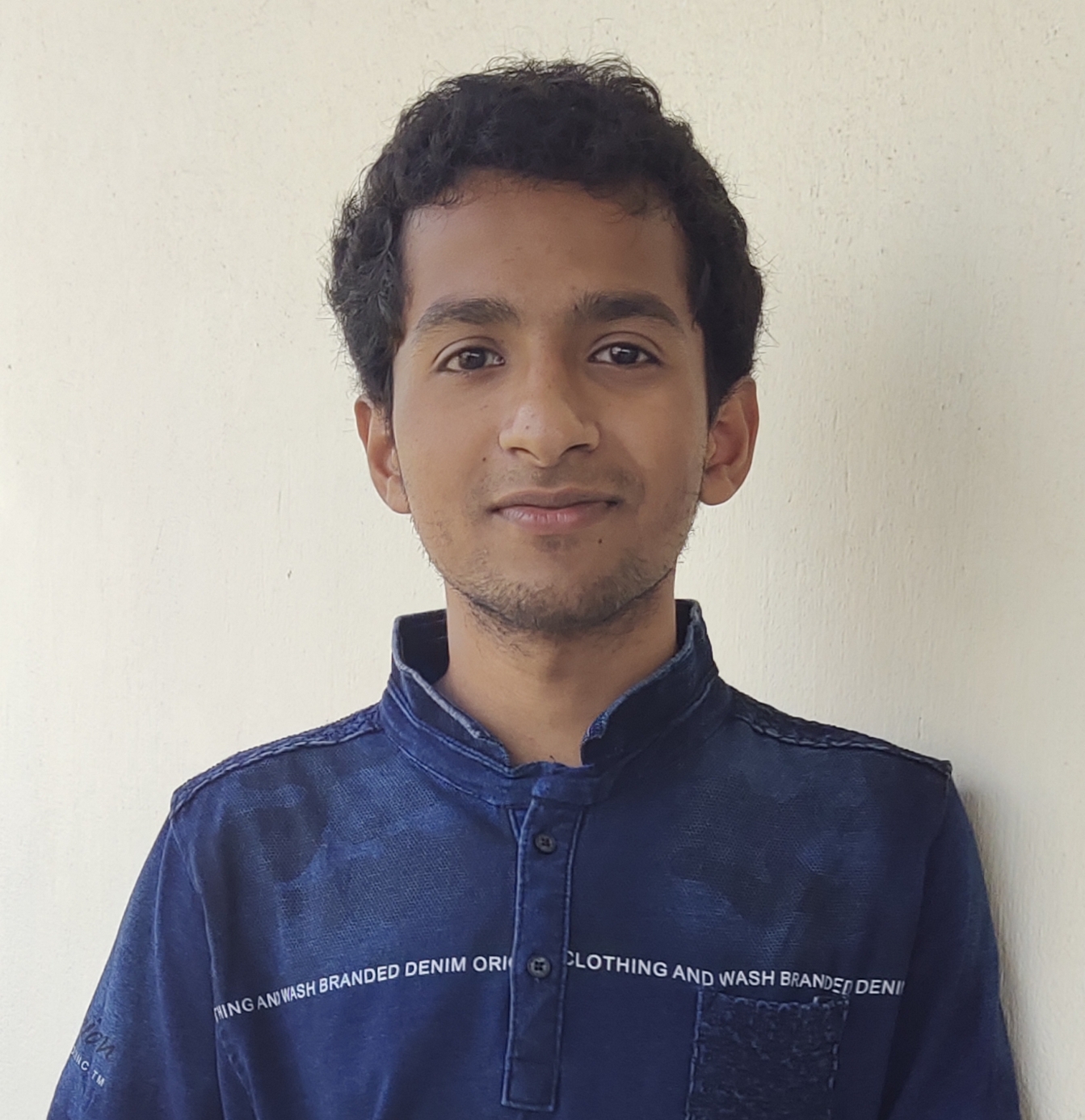}{The Author, currently is a student of BSc. Physics Hons. at Sripatsingh College, Under Kalyani University (India).}

\end{document}